\documentclass[12pt]{amsart}
\begin{document}

\title{Some additive relations in the Pascal triangle}
\author{A. V. Stoyanovsky}
\email{stoyan@mccme.ru}
\begin{abstract}
We derive some, seemingly new, curious additive relations in the Pascal triangle.
They arise in summing up the numbers in the triangle along some vertical line up to some place.
\end{abstract}

\maketitle

There are a lot of known additive relations among ($q$-)binomial coefficients, see, for
example, [1]. In this paper we derive a seemingly new curious relation obtained in a simplest
way, namely, by summation of the numbers of the Pascal triangle along a vertical line up to
some place.

It is well known that the sum of the numbers of the Pascal triangle,
$$
\begin{array}{ccccccccccc}
&&&&&1\\
&&&&1&&1\\
&&&1&&2&&1\\
&&1&&3&&3&&1\\
&1&&4&&6&&4&&1\\
1&&5&&10&&10&&5&&1\\
&&&&&\ldots
\end{array}
$$
along the $n$-th horizontal line equals $2^n$ (see, for example, [2]),
\begin{equation}
{n\choose 0}+{n\choose 1}+{n\choose 2}+\ldots+{n\choose n}=2^n.
\end{equation}
It is also known that the sum along the $k$-th diagonal up to some place equals the next
number in the next diagonal,
\begin{equation}
{k\choose k}+{k+1\choose k}+{k+2\choose k}+\ldots+{n\choose k}={n+1\choose k+1}.
\end{equation}
Another well known fact is that the sum along the $n$-th diagonal with the slope $1/3$
equals the $n$-th Fibonacci number,
\begin{equation}
{n\choose 0}+{n-1\choose 1}+{n-2\choose 2}+{n-3\choose 3}+\ldots=u_n,
\end{equation}
where
\begin{equation}
u_{n+2}=u_{n+1}+u_n,\ \ u_0=u_1=1
\end{equation}
are the Fibonacci numbers.

The main purpose of this note is to derive a formula for the sum along a vertical line up
to some place,
\begin{equation}
{n\choose k}+{n-2\choose k-1}+{n-4\choose k-2}+\ldots.
\end{equation}
This formula is the following.
\medskip

{\bf Theorem.} The sum (5) equals the alternated sum along the next diagonal with the slope
$1/3$, starting from the closest number, plus possibly $\pm 1$ depending on the vertical:
\begin{equation}
\begin{aligned}{}
&{n\choose k}+{n-2\choose k-1}+{n-4\choose k-2}+\ldots\\
&={n+1\choose k+1}-{n\choose k+2}+{n-1\choose k+3}-\ldots\\
&+\left\{
\begin{array}{l}
0,\  n-2k\le 0,\\
0,\  n-2k>0, n-2k=6p,6p+3,\\
-1,\  n-2k>0, n-2k=6p+1,6p+2,\\
1,\ n-2k>0, n-2k=6p-1,6p-2.
\end{array}\right.
\end{aligned}
\end{equation}

This theorem can be proved without big difficulties if one notes that both sides of
equality (6) satisfy the same recurrence relation as the numbers in the Pascal triangle,
\begin{equation}
\begin{aligned}{}
LHS(n,k)&=LHS(n-1,k)+LHS(n-1,k-1),\\
RHS(n,k)&=RHS(n-1,k)+RHS(n-1,k-1),
\end{aligned}
\end{equation}
and the same initial conditions (for $k=0$ and $k=n$).
More precisely, equalities (7) hold everywhere
outside the vertical line $n-2k=0$. On this line one should add $1$ to the right hand sides
of both equalities.

\end{document}